\documentclass[12pt, twoside]{article}
\usepackage{amsmath,amsthm,amssymb}
\usepackage{times}
\usepackage{enumerate}

\usepackage{epsfig}
\usepackage{subfigure}

\theoremstyle{definition}
\newtheorem{thm}{Theorem}[section]

\newtheorem{lem}[thm]{Lemma}
\newtheorem{defin}[thm]{Definition}
\newtheorem{rem}[thm]{Remark}

\newtheorem{mainthm}[thm]{Main Theorem}

\newtheorem*{xrem}{Remark}

\newcommand{\bbR}{{\mathbb R}}
\newcommand{\bbN}{{\mathbb N}}
\newcommand{\bbC}{{\mathbb C}}

\newcommand{\al}{\alpha}
\newcommand{\la}{\lambda}
\newcommand{\p}{\partial}
\newcommand{\be}{\beta}
\newcommand{\G}{\Gamma}

\numberwithin{equation}{section}

\newcommand{\subjclass}[1]{\bigskip\noindent\emph{2010 Mathematics Subject Classification:}\enspace#1}
\newcommand{\keywords}[1]{\noindent\emph{Keywords:}\enspace#1}

\begin{document}


\baselineskip=17pt


\title{A Generalized Neumann Solution for the Two-Phase Fractional Lam\'{e}-Clapeyron-Stefan Problem}

\author{ Sabrina D. ROSCANI\\
CONICET - Depto. Matem\'{a}tica,  FCEIA, Univ. Nac. de Rosario,\\
 Pellegrini 250, S2000BTP Rosario, Argentina  \\
sabrina@fceia.unr.edu.ar\\
Domingo A. TARZIA\\
  CONICET - Depto. Matem\'atica,
FCE, Univ. Austral,\\
 Paraguay 1950, S2000FZF Rosario, Argentina \\
dtarzia@austral.edu.ar}

\date{}

\maketitle


\begin{abstract}
We obtain a generalized Neumann solution for the  two-phase fractional Lam\'{e}-Clapeyron-Stefan problem for a semi-infinite material with constant boundary and initial conditions. In this problem, the two governing equations and a governing condition for the free boundary include a fractional time derivative in the Caputo sense of order $0<\al\leq 1$. When  $ \al \nearrow $ 1 we recover the classical Neumann solution for the two-phase Lam\'{e}-Clapeyron-Stefan problem given through the error function.


\end{abstract}

\section{Introduction}



The fractional diffusion equation has been treated by a number of authors (see  \cite{FM-Anal prop and apl of the W func, FM-The fundamental solution, Kilbas, Luchko3, Podlubny}) and, among the several applications that
have been studied, Mainardi \cite{FM-libro} studied the application
to the theory of linear viscoelasticity.\\

The free boundary problems for the one-dimensional heat equation  are problems linked to the processes
 of melting and freezing which have a latent heat-type condition at the interface connecting the velocity
 of the free boundary and the heat flux of the temperatures in both phases. This kind of problems have been widely studied
 (see \cite{Alexiades, Cannon, Car, Crank, Elliott, Gupta, Lame, Lunardini, Rubi, Stefan, Tar1, Tar2}).
\smallskip
In this paper, we deal with a two-phase Lam\'{e}-Clapeyron-Stefan problem for the time
fractional diffusion equation, obtained from the standard diffusion
equation by replacing the first order time-derivative by a
fractional derivative of order $\alpha \in (0,1) $ in the Caputo sense.\\

 We use here the definition introduced by Caputo in 1967 \cite{Caputo}, and we will call it \textsl {fractional derivative in the Caputo sense}, which is defined by

$$\,_{a} D^{\alpha}f(t)=\begin{cases} \frac{1}{\Gamma(n-\alpha)}\int^{t}_{a}(t-\tau)^{n-\alpha-1} f^{(n)}(\tau)d\tau & n-1<\al<n\\
f^{(n)}(t) & \al=n \end{cases}$$
where $\al >0$ is the order of derivation, $n \in \bbN$, $\G$ is the Gamma function defined by $\G(x)=\int_0^\infty  t^{x-1}e^{-t}dt$ and  $f $ is a differentiable function up to order $n$ in $\left[a,b\right]$.\\

 An interesting physical meaning of the fractional Stefan problems  is  discussed in \cite{Garra} and many authors were recently studying this kind of problems
 (see. \cite{Atkinson, Liu-Xu, RoSa1, RoSa2, Voller}). Some applications are linked to
 the behaviour in simulations of gas in polymer glasses (\cite{Gusev}) or propagation in porous media (\cite{Fellah}). In \cite{KZF2003} the classical Lam\'{e}-Clapeyron-Stefan problem was studied by using the fractional derivative of order $1/2$.\\

In this paper we consider the following two-phase fractional Lam\'{e}-Clapeyron-Stefan Problem

\begin{equation}\label{St-2fases}
\left\{\begin{array}{lll}
       (i)\hspace{0.5 cm}  _0D^{\al} u_2(x,t)=\lambda_2^2\dfrac{\partial^2u_2 }{\partial x^2}(x,t) &   0<x<s(t), \,  t>0, \,  0<\al<1, \\
        (ii)\hspace{0.4 cm} _0D^{\al} u_1(x,t)=\lambda_1^2\dfrac{\partial^2u_1 }{\partial x^2}(x,t) &   s(t)<x<\infty, \,  t>0, \,  0<\al<1, \\
			(iii)\hspace{0.3 cm}k_1u_{1x}(s(t),t)-k_2u_{2x}(s(t),t)=\rho l \,_0D^{\alpha}s(t) & t>0,\\
			(iv)\hspace{0.35 cm}	u_1(s(t),t)=u_2(s(t),t)=u_m &  t>0,  \\
			(v)\hspace{0.5 cm}	u_1(x,0)=u_1(+\infty,t)=u_i &   0<x<\infty, \\
       (vi)\hspace{0.35 cm} u_2(0,t)=u_0 & t>0,\\
			(vii)\hspace{0.2 cm}	 s(0)=0\end{array}\right.
\end{equation}

where $u_i<u_m<u_0$ and $\lambda_j^2=\dfrac{k_j}{\rho c_j}$, $j=1$(solid phase), $2$ (liquid phase).\\

In this problem, the two governing diffusion equations (\ref{St-2fases}-$ii$) and (\ref{St-2fases}-$i$) for $u_1$ and $u_2$ respectively, and the governing condition on the free boundary $s(t)$ (\ref{St-2fases}-$iii$) include a fractional time derivative in the Caputo sense of order $0<\al\leq 1$.
The goal of this paper is to obtain an explicit solution of the free boundary problem (\ref{St-2fases}), called a generalized Neumann solution with respect to the classical one given in \cite{Car}, \cite{Tarzia}, \cite{Weber}. This explicit solution is obtained through the Wright and Mainardi functions (\cite{FM-the M-Wright}). In Section $2$ a summary of some properties related to these special functions are given which will be useful in the next section. In Section $3$ the existence of a generalized Neumann solution is given and an open problem for the uniqueness is posed. Moreover, the classical Neumann solution for the two-phase Lam\'{e}-Clapeyron-Stefan problem for a semi-infinite material is well recovered by considering the limit when $\al \nearrow 1$.
\section{The Special Functions Involved}

\begin{defin} For every $z\in \mathbb{C}$ , $\al>-1$ and $\be\in \bbR$ the Wright function is defined by
\begin{equation}\label{W} W(z;\al;\be)=\sum^{\infty}_{k=0}\frac{z^{k}}{k!\G(\al k+\be)}.\end{equation}
\end{defin}
 This function will play a fundamental role in this paper. It is known that the Wright function is an entire function if $\al>-1$.\\

  Taking $\al=-\frac{1}{2}$ and $\be=\frac{1}{2}$, we get
 \begin{equation}  W\left(-z,-\frac{1}{2},\frac{1}{2}\right)=M_{1/2}(z)=\frac{1}{\sqrt{\pi}}e^{-z^2/4}.\end{equation}
where $M_{1/2}(z)$ is the Mainardi function (see \cite{FM-Anal prop and apl of the W func}), defined by
\begin{equation}\label{M} M_\nu (z)= W(-z,-\nu,1-\nu)=\sum^{\infty}_{n=0}\frac{(-z)^n}{n! \G\left( -\nu n+ 1-\nu \right)}, \quad z\in \bbC, \, \nu<1.  \end{equation}
which is a particular case of the Wright function. \\

 Due to the uniform convergence of the series on compact sets, we have (see \cite{Wright 2})
 \begin{equation}\label{derivada de W} \frac{\p}{\p z} W(z,\al,\be) = W(z,\al,\al+\be).  \end{equation}
 Then, for $x\in \bbR^+_0$, and taking account that
 \begin{equation}\label{W(-inf,-al 2,1)} W(-\infty,-\frac{\al}{2},1)=0, \qquad  \text{ if } \
  \al\in (0,2),
  \end{equation}
 we have
 $$ W\left(-x,-\frac{1}{2},1\right)=  W\left(-x,-\frac{1}{2},1\right)-
 W\left(-\infty,-\frac{1}{2},1\right) = \int_{\infty}^x \left( \frac{\p}{\p x}
W\left(-\xi,-\frac{1}{2},1\right)\right)d\xi= $$
$$\hspace{3 cm}= \int_{\infty}^x-W\left(-\xi,-\frac{1}{2},\frac{1}{2}\right)d\xi = \int_x^{\infty}
W\left(-\xi,-\frac{1}{2},\frac{1}{2}\right)d\xi=\int_x^{\infty}
\frac{1}{\sqrt{\pi}}e^{-\xi^2/4}d\xi=$$
 $$\hspace{-4.6 cm} = \frac{2}{\sqrt{\pi}}\int_{x/2}^{\infty} e^{-\xi^2}d\xi
 = \mbox{erfc\,}\left(\frac{x}{2}\right), $$
 that is,
\begin{equation}  W\left(-x,-\frac{1}{2},1\right)= \mbox{erfc\,}\left(\frac{x}{2}\right), \quad
1-W\left(-x,-\frac{1}{2},1\right)= \mbox{erf\,}\left(\frac{x}{2}\right).\end{equation}

\noindent where {\it erf} and {\it erfc} are the error and complementary error functions.\\

\noindent The next two propositions were proved in \cite{RoSa1}.

\begin{lem}\label{M pos y decrec} If  $\,  0<\al<1$, then:
\begin{enumerate}
\item $M_{\al/2}(x) $ is a positive and strictly decreasing positive function in $\bbR^+$  such that $M_{\al/2}(x)<\frac{1}{\G\left(1-\frac{\al}{2}\right)}$;
\item  $W\left(-x,-\frac{\al}{2},1\right)$ is a positive and strictly decreasing function in  $\bbR^+$ such that $0<W\left(-x,-\frac{\al}{2},1\right)\leq 1, $ $\, \forall \, x \in \bbR^+_0$.
\end{enumerate}
\end{lem}

\begin{lem}\label{conv M al/2 cuando al tiende a 1}If $x\in \bbR^+_0$ and $\al \in (0,1)$ then:
\begin{enumerate}

\item $\lim\limits_{\al\nearrow 1}M_{\al/2}\left(x\right)=M_{1/2}(x)=\frac{e^{-\frac{x^2}{4}}}{\sqrt{\pi}}$;
\item $\lim\limits_{\al\nearrow 1}\left[1-W\left(-x,-\frac{\al}{2},1\right)\right]=\frac{1}{\sqrt{\pi}} erf\left(\frac{x}{2}\right) .$

\end{enumerate}

\end{lem}

Due to the results in \cite{Wright}, the following assertions are true

\begin{equation}\label{lim W(-x,-al/2,1)=0 y M=0}
\lim_{x \rightarrow \infty}W\left(-x,-\frac{\al}{2},1\right)=0
\quad \text{ and } \quad \lim_{x \rightarrow \infty}M_{\al/2}\left(x\right)=0. \quad
\end{equation}

Let us work on some problems in the first quadrant. It is known that (see \cite{FM-The fundamental solution})
\begin{equation}{\label{sol gral para DFE}}
u(x,t)=\int^{\infty}_{-\infty} \frac{t^{-\frac{\al}{2}}}{2\lambda}
M_{\frac{\al}{2}}\left(\left|x-\xi\right|\lambda^{-1}t^{-\frac{\al}{2}}\right) f(\xi)d\xi
\end{equation}
 is a solution for the fractional diffusion problem
\begin{equation}{\label{P1}}
\left\{\begin{array}{lll}
          _{0}D^{\al} u(x,t)=\lambda^2\dfrac{\partial^2u }
          {\partial x^2}(x,t) &   -\infty<x<\infty, \,  t>0, \,  0<\al<1, \\
         u(x,0)=f(x) &  -\infty<x<\infty.
      \end{array}\right.\end{equation}
Using this fact, it is easy to see that
\begin{equation}{\label{sol DFE 1er cuad}}
 v(x,t)=\frac{1}{2\lambda t^{\frac{\al}{2}}}\int^{\infty}_{0}\left[ M_{\frac{\al}{2}}\left(\frac{|x-\xi |}{
\lambda t^{\frac{\al}{2}}}\right)-M_{\frac{\al}{2}}\left(\frac{x+\xi}{
\lambda t^{\frac{\al}{2}}}\right) \right]
 f_0 \, d\xi
\end{equation}
is a solution for the fractional diffusion problem
\begin{equation}{\label{P1}}
\left\{\begin{array}{lll}
          _{0}D^{\al} v(x,t)=\lambda^2\dfrac{\partial^2 v }{\partial x^2}(x,t) &   0<x<\infty, \,  t>0, \,  0<\al<1, \\
         v(x,0)=f_0 &   0<x<\infty, \\
         v(0,t)=0 & t>0.                                  \end{array}\right.\end{equation}

\smallskip
\noindent
An equivalent expression of (\ref{sol DFE 1er cuad}) is given by
$$\hspace{-7 cm} v(x,t)=\frac{1}{2\lambda t^{\frac{\al}{2}}}\int^{\infty}_{0}\left[ M_{\frac{\al}{2}}\left(\frac{|x-\xi |}{
\lambda t^{\frac{\al}{2}}}\right)-M_{\frac{\al}{2}}\left(\frac{x+\xi}{
\lambda t^{\frac{\al}{2}}}\right) \right]
 f_0 d\xi $$
$$\hspace{0.5 cm}=\frac{f_0}{2}\left[\int^{x}_{0}\frac{1}{\lambda t^{\frac{\al}{2}}}M_{\frac{\al}{2}}\left(\frac{x-\xi }{
\lambda
t^{\frac{\al}{2}}}\right)d\xi+\int^{\infty}_{x}\frac{1}{\lambda
t^{\frac{\al}{2}}}M_{\frac{\al}{2}}\left(\frac{\xi-x }{ \lambda
t^{\frac{\al}{2}}}\right)d\xi\right.\left.-\int^{\infty}_{0}\frac{1}{\lambda
t^{\frac{\al}{2}}}M_{\frac{\al}{2}}\left(\frac{x+\xi}{ \lambda
t^{\frac{\al}{2}}}\right)d\xi\right]$$
$$\hspace{0.2 cm}=\frac{f_0}{2}\left[-W\left(-\frac{x}{\lambda t^{\frac{\al}{2}}},-\frac{\al}{2},1\right)+2 -W\left(-\frac{x}{\lambda
 t^{\frac{\al}{2}}},-\frac{\al}{2},1\right)\right]= f_0\left[1-W\left(-\frac{x}{\lambda
t^{\frac{\al}{2}}},-\frac{\al}{2},1\right)\right].
 $$
Analogously we can check that
\begin{equation}\label{sol W}
w(x,t)=g_0 \, W\left(-\frac{x}{\la t^{\al/2}},-\frac{\al}{2},1\right)
\end{equation}
is a solution for the fractional diffusion problem
\begin{equation}{\label{P1}}
\left\{\begin{array}{lll}
         _{0}D^{\al} w(x,t)=\lambda^2\dfrac{\partial^2 w }{\partial x^2}(x,t) &   0<x<\infty, \,  t>0, \,  0<\al<1, \\
         w(x,0)=0 &   0<x<\infty, \\
         w(0,t)=g_0 & t>0.                                  \end{array}\right.\end{equation}

\medskip

\section{The Two-Phase Fractional Lam\'{e}-Clapeyron-Stefan Problem}

\noindent Hereinafter we will call $D^\al$ to the fractional derivative in the Caputo sense of extreme $a=0$, $_0D^\al$.\\

\noindent Let us return to problem (\ref{St-2fases}). Taking into account the previous section and the method developed  in \cite{RoSa1}, the following explicit solution is obtained.

\begin{thm}
An explicit solution for the two-phase Lam\'{e}-Clapeyron-Stefan problem (\ref{St-2fases}) is given by

\begin{equation}\label{sol st-2f}
\left\{\begin{array}{l} u_2(x,t)=u_0-(u_0-u_m)\frac{1-W\left(-\frac{x}{\lambda_2 t^{\al/2}},-\frac{\al}{2},1\right)}{1-W\left(-\xi \lambda,-\frac{\al}{2},1\right)}\\
        u_1(x,t)=u_i+(u_m-u_i)\frac{W\left(-\frac{x}{\lambda_1 t^{\al/2}},-\frac{\al}{2},1\right)}{W\left(-\xi,-\frac{\al}{2},1\right)} \\
				 s(t)=\xi  \lambda_1 t^{\al/2} \end{array}\right.
\end{equation}
 where $\xi$ is a solution to the equation
\begin{equation}\label{eq xi}
F(x)=\frac{\G(1+\frac{\al}{2})}{\G(1-\frac{\al}{2})}x, \, x>0   \end{equation}

\noindent and the function $F:\bbR^+_0\rightarrow \bbR$ is defined by
\begin{equation}\label{F}
F(x)=\frac{k_2(u_0-u_m)}{\rho l \lambda_1 \lambda_2}F_1(\lambda x)- \frac{k_1(u_m-u_i)}{\rho l \lambda_1^2}F_2(x)
\end{equation}

\noindent with
\begin{equation}\label{F1 y F2}  F_1(x)=\frac{M_{\al/2}(x)}{1-W\left(-x,-\frac{\al}{2},1\right)} \,, \,\,
F_2(x)=\frac{M_{\al/2}(x)}{W\left(-x,-\frac{\al}{2},1\right)},\,\, \lambda=\frac{\lambda_1}{\lambda_2}>0. \end{equation}

\end{thm}

\proof
The following solution is proposed
\begin{equation}\label{Sol-propuesta}
\left\{\begin{array}{l} u_2(x,t)=A+B\left[1-W\left(-\frac{x}{\lambda_2 t^{\al/2}},-\frac{\al}{2},1\right)\right]\\
        u_1(x,t)=C+D\left[1-W\left(-\frac{x}{\lambda_1 t^{\al/2}},-\frac{\al}{2},1\right)\right] \\
				 s(t)=\xi  \lambda_1 t^{\al/2}\end{array}\right.
\end{equation}

where $A, B, C, D$ and $\xi >0$ must be determined. \\
According with the results in the previous section and the linearity of the fractional derivative $D^\al$, functions $u_2$ and $u_1$ are solutions of the fractional diffusion equations (\ref{St-2fases}-i) and (\ref{St-2fases}-ii) respectively.

From conditions (\ref{St-2fases}-$iv$) and (\ref{St-2fases}-$vi$)    we have,
\begin{equation}\label{St-2f-1}
u_2(0,t)=A+B\left[1-W\left(0,-\frac{\al}{2},1\right)\right]=u_0
\end{equation}

\begin{equation}\label{St-2f-2}
u_2(s(t),t)=u_0+B\left[1-W\left(-\xi \frac{\lambda_1}{\lambda_2},-\frac{\al}{2},1\right)\right]=u_m.
\end{equation}
and therefore we obtain:
\begin{equation}\label{A y B} A=u_0 , \quad \text{and} \quad  B=-\frac{u_0-u_m}{1-W\left(-\xi \lambda,-\frac{\al}{2},1\right)}.\end{equation}
So,

\begin{equation}\label{u2a}
u_2(x,t)=u_0-(u_0-u_m)\frac{1-W\left(-\frac{x}{\lambda_2 t^{\al/2}},-\frac{\al}{2},1\right)}{1-W\left(-\xi \lambda,-\frac{\al}{2},1\right)}<u_0,
\end{equation}

or equivalently

\begin{equation}\label{u2b}
u_2(x,t)=u_m+(u_0-u_m)\frac{W\left(-\frac{x}{\lambda_2 t^{\al/2}},-\frac{\al}{2},1\right)-W\left(-\xi \lambda,-\frac{\al}{2},1\right)}{1-W\left(-\xi \lambda,-\frac{\al}{2},1\right)}.
\end{equation}

Taking into account the results in Proposition \ref{M pos y decrec}, (\ref{u2a}) and (\ref{u2b}) it is easy to see that
\begin{equation}\label{u_2>u_m} u_m<u_2(x,t)<u_0, \quad 0<x<s(t), \, t>0 .
\end{equation}

From conditions (\ref{St-2fases}-$v$) and (\ref{St-2fases}-$iv$)    we have,

\begin{equation}\label{St-2f-3}
u_1(x,0)=C+D\left[1-W\left(-\infty,-\frac{\al}{2},1\right)\right]=C+D=u_i,
\end{equation}
\begin{equation}\label{St-2f-4}
u_1(s(t),t)=C+D\left[1-W\left(-\xi,-\frac{\al}{2},1\right)\right]=u_m,
\end{equation}
and therefore we get:

\begin{equation}\label{St-2f-5}
C=u_i+\frac{u_m-u_i}{W\left(-\xi,-\frac{\al}{2},1\right)}, \quad D=-\frac{u_m-u_i}{W\left(-\xi,-\frac{\al}{2},1\right)}.
\end{equation}

Accordingly,

\begin{equation}\label{u1a}
u_1(x,t)=u_i+(u_m-u_i)\frac{W\left(-\frac{x}{\lambda_1 t^{\al/2}},-\frac{\al}{2},1\right)}{W\left(-\xi,-\frac{\al}{2},1\right)},
\end{equation}
or equivalently

\begin{equation}\label{u1b}
u_1(x,t)=u_m-(u_m-u_i)\left[1-\frac{W\left(-\frac{x}{\lambda_1 t^{\al/2}},-\frac{\al}{2},1\right)}{W\left(-\xi,-\frac{\al}{2},1\right)}\right].
\end{equation}

Taking into account  Proposition \ref{M pos y decrec}, (\ref{u1a}) and (\ref{u1b}) we obtain

\begin{equation}\label{u_i<u1<u_m}
u_i<u_1(x,t)<u_m, \quad x>s(t)=\xi\lambda_1 t^{\al/2}, \, t>0.
\end{equation}

In order to determine  $\xi>0$, let us work with the ``fractional Lam\'{e}-Clapeyronn-Stefan condition'' (\ref{St-2fases}-$iii$). From (\ref{W}) and (\ref{derivada de W}) we have

$$  u_{2x}(x,t)=\frac{B}{\lambda_2 t^{\al/2}}M_
{\al/2}\left(\frac{x}{\lambda_2 t^{\al/2}}\right) ,\quad  u_{1x}(x,t)=\frac{D}{\lambda_1 t^{\al/2}}M_
{\al/2}\left(\frac{x}{\lambda_1 t^{\al/2}}\right), $$

which evaluated on $(s(t),t)$, gives

\begin{equation}\label{cond st 1}
u_{2x}(s(t),t)=\frac{B}{\lambda_2 t^{\al/2}}M_
{\al/2}\left(\lambda \xi\right) ,\quad u_{1x}(s(t),t)=\frac{D}{\lambda_1 t^{\al/2}}M_
{\al/2}\left(\xi\right).
\end{equation}

Taking into account that (\cite{Podlubny})
$$ D^{\al}(t^{\be})=\frac{\G(\be+1)}{\G(1+\be-\al)}t^{\be-\al} \quad \text{if } \be>-1, $$
it results that
\begin{equation}\label{D^al s}
D^{\al}s(t)=D^{\al}( \xi \lambda_1 t^{\al/2})=\lambda_1\xi D^{\al}(t^{\al/2})=\lambda_1 \xi \frac{\G(1+\frac{\al}{2})}{\G(1-\frac{\al}{2})}t^{-\al/2}.
\end{equation}

Replacing (\ref{cond st 1}) and (\ref{D^al s}) in the fractional condition (\ref{St-2fases}-iii), we get for the unknown coefficient $\xi >0 $ the following equation:

$$ 	 k_1u_{1x}(s(t),t)-k_2u_{2x}(s(t),t)=\rho l \,D^{\alpha}s(t) \Leftrightarrow  $$
$$ -k_1 \frac{u_m-u_i}{1-W\left(- \xi, -\frac{\al}{2},1\right)} \frac{1}{\lambda_1 t^{\al/2}}M_{\al/2}(\xi)+k_2 \frac{u_0-u_m}{1-W\left(-\lambda \xi, -\frac{\al}{2},1\right)}\frac{1}{\lambda_2 t^{\al/2}}M_{\al/2}(\lambda \xi)= \rho l \lambda_1 \xi \frac{\G(1+\frac{\al}{2})}{\G(1-\frac{\al}{2})}t^{-\al/2} \Leftrightarrow$$
$$ \frac{k_2(u_0-u_m)}{\lambda_2}\frac{M_{\al/2}(\lambda \xi)}{1-W\left(-\lambda \xi , -\frac{\al}{2},1\right)}-
 \frac{k_1(u_m-u_i)}{\lambda_1}\frac{M_{\al/2}( \xi)}
 {W\left(-\xi ,-\frac{\al}{2},1\right)}=\xi \rho l \lambda_1 \frac{\G(1+\frac{\al}{2})}{\G(1-\frac{\al}{2})} \Leftrightarrow $$
$$ \frac{k_2(u_0-u_m)}{\rho l \lambda_1 \lambda_2}F_1(\lambda \xi)- \frac{k_1(u_m-u_i)}{\rho l \lambda_1^2}F_2(\xi)=
\frac{\G(1+\frac{\al}{2})}{\G(1-\frac{\al}{2})}\xi \Leftrightarrow $$
\begin{equation}  \Leftrightarrow  F(\xi)=\frac{\G(1+\frac{\al}{2})}{\G(1-\frac{\al}{2})}\xi\, ; \end{equation}
that is, the equation (\ref{eq xi}) holds, where $F$, $F_1$ and $F_2$ where defined in (\ref{F}) and (\ref{F1 y F2})  respectively.\\

In order to guarantee the existence of a solution of the equation (\ref{eq xi}), we will study the behavior of the functions $F$, $F_1$ and $F_2$. From Proposition \ref{M pos y decrec} and (\ref{lim W(-x,-al/2,1)=0 y M=0}), it results that
\begin{equation} \label{caract F_1} F_1  \text{ is a positive decreasing function, } \, F_1(0^+)=\infty, \, \text{ and } F_1(+\infty)=0 \end{equation}
and
\begin{equation}\label{caract F_2} F_2 \, \text{is a positive function and } F_2(0)=\frac{1}{\G{(1-\al/2)}}.\end{equation}

Let us prove that \begin{equation}\label{F_2(+infty)=+infty}F_2(+\infty)=+\infty.\end{equation}

 In \cite{Wright 3, Wright} the asymptotic expansion for $x\rightarrow \infty$ of the Wright function was studied, and an interesting summary of these results can be founded in  \cite{Wong}, from where we can say that if $\al \in (0,1)$ we have
$$ M_{\al/2}(x)= \left(\frac{\al}{2} x\right)^{-\frac{1-\al}{2-\al}}\exp \left\{ \left(1-\frac{2}{\al}\right)\left(\frac{\al}{2} x\right)^{\frac{1}{1-\al/2}}\right\}\left[a_0+\mathcal{O}\left(\left(\frac{\al}{2} x\right)^{-\frac{1}{1-\al/2}}\right)\right], \quad a_0=\frac{1}{\sqrt{2\pi(1-\al/2)}} $$

Therefore
\begin{equation}\label{comp asint M}
M_{\al/2}(x)\sim
b(\al)x^{-\frac{1-\al}{2-\al}}\exp \left\{ -c(\al)x^{\frac{1}{1-\al/2}}\right\}
\end{equation}
where  $b(\al)=\frac{1}{\sqrt{2\pi(1-\al/2)}}\left(\frac{\al}{2}\right)^{-\frac{1-\al}{2-\al}}>0$ and $c(\al)=\frac{2-\al}{2}\left(\frac{\al}{2} \right)^{\frac{1}{1-\al/2}}>0$.\\

On the other hand
$$
W\left(-x,-\frac{-\al}{2},1\right)= \left(\frac{\al}{2} x\right)^{-\frac{1}{2-\al}}\exp \left\{ \left(1-\frac{2}{\al}\right)\left(\frac{\al}{2} x\right)^{\frac{1}{1-\al/2}}\right\}\left[a_0+\mathcal{O}\left(\left(\frac{\al}{2} x\right)^{-\frac{1}{1-\al/2}}\right)\right] ,$$

therefore
\begin{equation}\label{comp asint W}
W\left(-x,-\frac{-\al}{2},1\right)\sim
d(\al)x^{-\frac{1}{2-\al}}\exp \left\{ -c(\al)x^{\frac{1}{1-\al/2}}\right\}
\end{equation}

where  $d(\al)=\frac{1}{\sqrt{2\pi(1-\al/2)}}\left(\frac{\al}{2}\right)^{-\frac{1}{2-\al}}>0$.\\

From (\ref{comp asint M}) and (\ref{comp asint W}), we have
\begin{equation}\label{comp asint F_2}
F_2(x)\sim \left(\frac{\al}{2}\right)^{\frac{\al}{2-\al}}x^{\frac{\al}{2-\al}}, \quad \text{as } x\rightarrow \infty
\end{equation}
and then (\ref{F_2(+infty)=+infty}) holds.\\

Now, from Proposition \ref{M pos y decrec} and properties  (\ref{caract F_2}) and (\ref{F_2(+infty)=+infty}), we can ensure that

\begin{equation}\label{F(0) y F(infty)}
F \text{ is a continuous function, }\,  F(0^+)=+\infty \quad \text{ and } \quad F(+\infty)=-\infty.
\end{equation}

Therefore, there exists at least one $\xi>0$ which is solution of the equation (\ref{eq xi}). Finally, we are able to state that (\ref{sol st-2f}) is a solution to the free boundary problem (\ref{St-2fases}).

\endproof

\begin{rem} We will denote (\ref{sol st-2f})-(\ref{eq xi}) as  the generalized Neumann solution of the  two-phase fractional Lam\'{e}-Clapeyron-Stefan problem (\ref{St-2fases}).

\end{rem}

\begin{thm}The limit when $\al\nearrow 1$ of the generalized Neumann solution (\ref{sol st-2f})-(\ref{eq xi}) is the classical Neumann solution for the two-phase Lam\'{e}-Clapeyron-Stefan problem.

\end{thm}

\proof

We denote $u_1^\al$, $u_2^\al$ and $s_\al$ as the functions defined in (\ref{sol st-2f}), and $\xi_\al$ the solution of the equation (\ref{eq xi}) for each $0<\al<1$. Now, we analyze the convergence of (\ref{sol st-2f}) when  $\al \nearrow 1$.   Applying Proposition  \ref{conv M al/2 cuando al tiende a 1} we obtain

\begin{equation}\label{conv 1}\lim_{\al\nearrow 1}u_1^\al(x,t)=
     u_i+(u_m-u_i)\lim_{\al\nearrow 1}\frac{W\left(-\frac{x}{\lambda_1 t^{\al/2}},-\frac{\al}{2},1\right)}{W\left(-\xi,-\frac{\al}{2},1\right)}=u_i+(u_m-u_i)\frac{\mbox{erfc\,} \left(\frac{x}{2\lambda_1\sqrt{t}}\right)}{\mbox{erfc\,} \left(\frac{\xi}{2}\right)} \end{equation}

\begin{equation}\label{conv 2}\lim_{\al\nearrow 1}u_2^\al(x,t)=
     u_0-(u_0-u_m)\lim_{\al\nearrow 1}\frac{1-W\left(-\frac{x}{\lambda_2 t^{\al/2}},-\frac{\al}{2},1\right)}{1-W\left(-\xi\lambda,-\frac{\al}{2},1\right)}=u_0-(u_0-u_m)\frac{\mbox{erf\,}\left(\frac{x}{2\lambda_2\sqrt{t}}\right)}{\mbox{erf\,} \left(\frac{\xi\lambda}{2}\right)} \end{equation}

     \begin{equation}\label{conv 3 }\lim_{\al\nearrow 1} s_\al(t)=\lim_{\al\nearrow 1}\xi_\al  \lambda_1 t^{\al/2}=\xi_1\lambda_1 \sqrt{t}=2\mu \lambda_1 \sqrt{t}  \end{equation}

where $ \mu=\frac{\xi_1}{2}$  is a solution to the equation
  \begin{equation}\label{conv 4} \frac{k_2(u_0-u_m)}{\rho l \lambda_1 \lambda_2}\frac{\exp\left\{ -\lambda^2 \mu^2 \right\}}{\sqrt{\pi} \mbox{erf\,}\left( \lambda \mu \right)}- \frac{k_1(u_m-u_i)}{\rho l \lambda_1^2}\frac{\exp \left\{-\mu^2\right\}}{\sqrt{\pi}\mbox{erfc\,} \left(\mu\right)}=
\mu, \quad \mu>0. \end{equation}
The expressions (\ref{conv 1})-(\ref{conv 4}) give us the classical Neumann solution, given in \cite{Car, Tarzia, Weber}, for the two-phase Lam\'{e}-Clapeyron-Stefan problem defined by the following equations, and constant boundary and initial conditions:
\begin{equation}
\left\{\begin{array}{lll}
         \dfrac{\partial u_2 }{\partial t}(x,t)=\lambda_2^2\dfrac{\partial^2u_2 }{\partial x^2}(x,t) &   0<x<s(t), \,  t>0, \,  \\
         \dfrac{\partial u_1 }{\partial t}(x,t)=\lambda_1^2\dfrac{\partial^2u_1 }{\partial x^2}(x,t) &   s(t)<x<\infty, \,  t>0, \, \\
				 k_1u_{1x}(s(t),t)-k_2u_{2x}(s(t),t)=\rho l \,\dot{s}(t) & t>0.\\
				u_1(s(t),t)=u_2(s(t),t)=u_m &  t>0,  \\
				u_1(x,0)=u_i &   0<x<\infty \\
         u_2(0,t)=u_0 & t>0.\\
				 s(0)=0\end{array}\right.
\end{equation}

\endproof

\begin{rem} It is an open problem to prove that $F_2$ is an increasing function, which is a sufficient condition to could ensure the uniqueness of the solution to equation (\ref{eq xi}). By using Maple we show below some graphs for different values of $0<\al<1$, from which it can be seen that $F_2$ is an increasing function on $\bbR^+$.

\begin{figure}[th]
\centering
\mbox {\subfigure[$F_2$ is an increasing function for $\al=1/16, 1/8, 1/4, 3/8$ and $1/2.$ ]{\epsfysize=75mm \epsfbox
{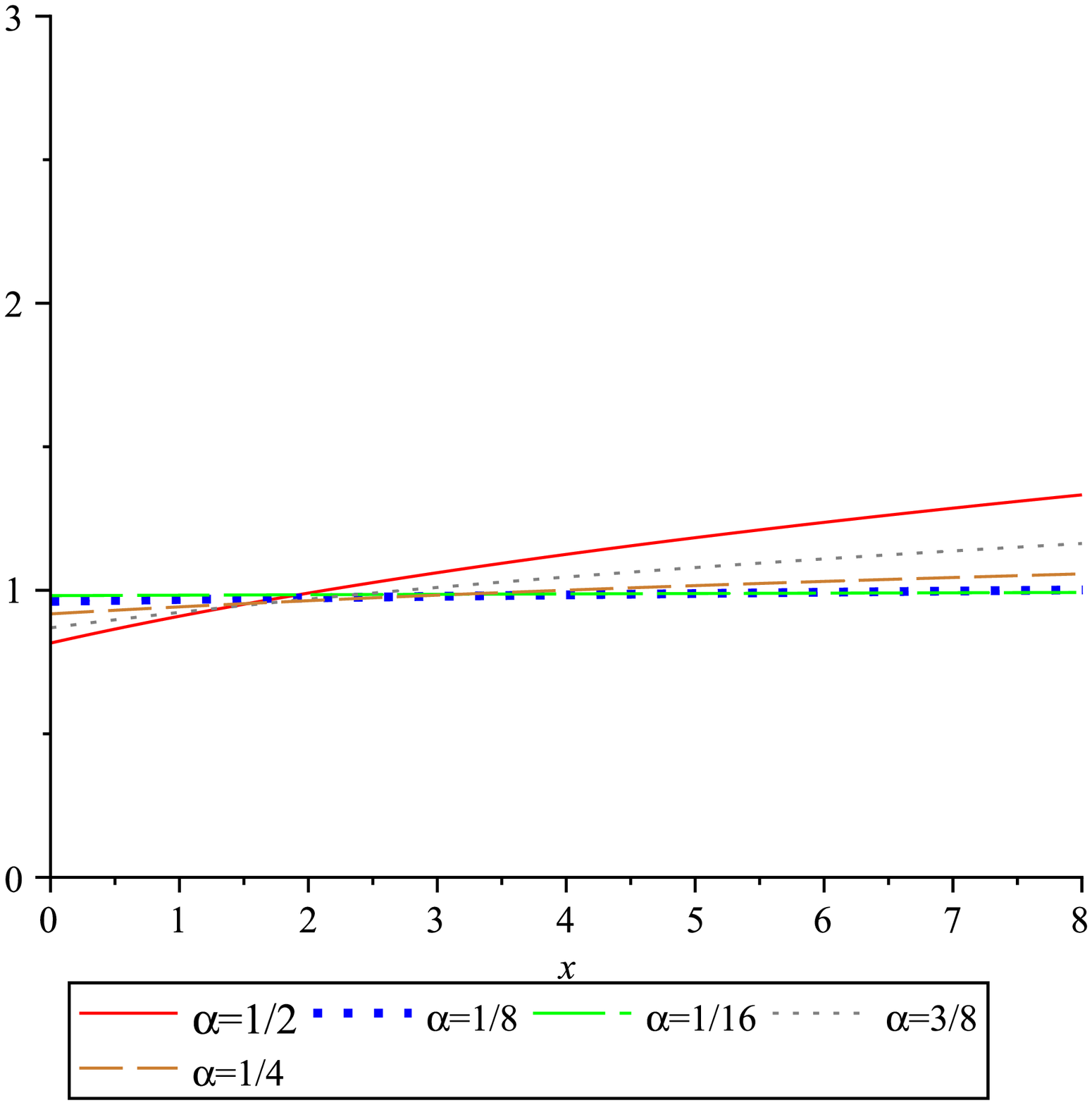}} \quad \quad \quad
\subfigure[$F_2$ is an increasing function for: $\al=1/2, 5/8, 7/8, 3/4$ and $15/16.$ ] {\epsfysize=75mm \epsfbox{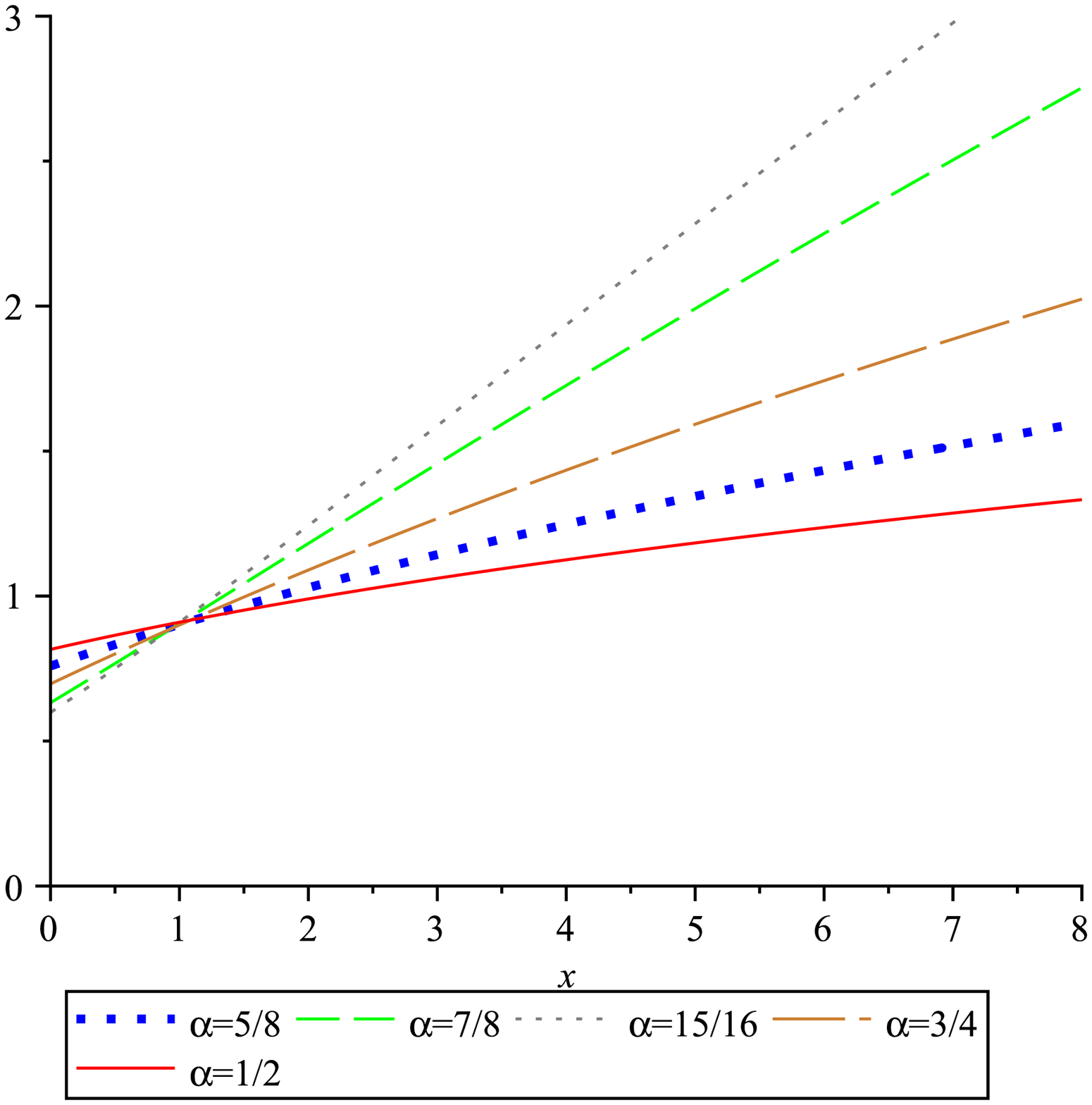}} }
\end{figure}


\end{rem}

\section{Conclusions}

By using the Wright and Mainardi functions and the fractional error function $1-W(-x,-\al/2,1)$, a generalized Neumann solution for the two-phase fractional Lam\'{e}-Clapeyron-Stefan problem is obtained for each $0<\al<1$. Moreover, the classical Neumann solution is recovered through the limit when $\al \nearrow 1$.

\subsection*{Acknowledgments}
\noindent This paper has been sponsored by the Projects ANPCyT PICTO AUSTRAL 2008 N. 73, PIP N. 0534 from CONICET-UA, Rosario, and ING349, from Universidad Nacional de Rosario, Argentina.

\end{document}